\documentstyle[amsmath]{article}
\begin{document}
\input{amssym}
\def\be{\begin{eqnarray*}}
\def\ee{\end{eqnarray*}}
\def\di{\displaystyle}
\title{Web geometry of a system of $n$ first order autonomous ordinary differential equations}
\author{Mehdi Nadjafikhah\thanks{School of Mathematics, Iran University of Science and Technology, Narmak-16, Tehran, I.R. Iran. e-mail: m\_nadjafikhah@ius.ac.ir}}
\maketitle
\begin{abstract}
Let $dx_i/dt=f_i(x_1,\cdots,x_n)$, ($i=1,\cdots,n$) be a system of
$n$ first order autonomous ordinary differential equations. We use
E. Cartan's equivalence method to study the invariants of this
system under diffeomorphisms of the form
$\Phi(t,x_1,\cdots,x_n)=(\varphi_0(t),\varphi_1(x_1),\cdots,\varphi_1(x_1))$.
\end{abstract}
\section{Introduction}
The method of equivalence of E. Cartan (see \cite{[2]}, \cite{[3]}
and \cite{[6]}) provides a powerful tool for constructing
differential invariants which solve the problem of deciding when
two geometric objects are really the same up to some preassigned
group of coordinate transformations. The solution of this problem
in terms of differential invariants goes back to S. Lie, but the
essential contributions of Cartan was the construction of an
adapted coframe whose structure equations yield differential
invariants. The analysis of the structure equations provides
classification results and yields a natural way of giving
invariant characterizations of the special models. In \cite{[2.5]}
R.B. Gardner gave some examples of solving these problems. For
example, he give the local equivalence problem for $dy/dx=f(x,y)$
under diffeomorphisms of the form
$\Phi(x,y)=(\varphi(x),\psi(y))$. We generalize this problem to a system of
$n$ first order autonomous ordinary differential equations.

In this paper we present a solution to the local equivalence
problem for
\begin{eqnarray}\label{eq:1}
\frac{dx_i}{dt}=f_i(x_1,\cdots,x_n),\qquad i=1,\cdots,n
\end{eqnarray}
under the group of coordinate transformations defined by
\begin{eqnarray}\label{eq:2}
\Phi(t,x_1,\cdots,x_n)=(\varphi_0(t),\varphi_1(x_1),\cdots,\varphi_1(x_1))
\end{eqnarray}
This is called {\it the pseudo-group of web transformations}.
\section{Setting the problem}
Given $(U,t,x_1,\cdots,x_n)$ and $(V,T,X_1,\cdots,X_n)$ open sets
with coordinates in ${\Bbb R}^{n+1}$ and ordinary differential
equations (\ref{eq:1}) on $U$ and
\begin{eqnarray}\label{eq:3}
\frac{dX_i}{dT}=F_i(X_1,\cdots,X_n),\qquad i=1,\cdots,n
\end{eqnarray}
on $V$. The usual symmetries of the ordinary differential
equations are the diffeomorphisms $\Phi:{\Bbb R}^{n+1}\to{\Bbb
R}^{n+1}$ which map integral curves into integral curves, that is,
satisfy
\begin{eqnarray}\label{eq:4}
\Phi^*(dX_i-F_i\,dT)=\sum_{j=1}^na_{ij}(dx_i-f_i\,dt)
\end{eqnarray}
where $[a_{ij}]:{\Bbb R}^{n+1}\to{\rm GL}(n,{\Bbb R})$ is a smooth
function. We also have the condition $\Phi^*(dT)=a\,dt$, where
$a:{\Bbb R}^{n+1}\to{\Bbb R}^*$ is a smooth function. Thus, we
have the following condition on Jacobian matrix of diffeomorphism
$\Phi$:
\begin{eqnarray}\label{eq:6}
\Phi^*\Omega=g.\omega,
\end{eqnarray}
where $\Omega=(\Omega^0,\cdots,\Omega^n)^t$,
$\omega=(\omega^0,\cdots,\omega^n)^t$,
\begin{align}\label{eq:7}
\Omega^0=dT,\;\; \Omega^i=dX_i-F_i\,dT,\;\; \omega^0=dt,\;\;
\omega^i=dx_i-f_i\,dt,\;\;\;(i=1,\cdots,n)
\end{align}
and $g:{\Bbb R}^{n+1}\to G$ is a smooth function, with
\begin{eqnarray}\label{eq:8}
G=\left\{\left(\begin{array}{cc}a&0\\0&A\end{array}\right)\;|\;a\in{\Bbb
R}^*,\; A\in{\rm GL}(n,{\Bbb R})\right\}.
\end{eqnarray}
We also have the condition on the Jacobian of the diffeomorphism
that
\begin{eqnarray}\label{eq:9}
\Phi^*\Theta=h.\theta
\end{eqnarray}
where $\Theta=(\Theta^0,\cdots,\Theta^n)^t$,
$\theta=(\theta^0,\cdots,\theta^n)^t$,
\begin{align}\label{eq:10}
\Theta^0=dT,\;\; \Theta^i=dX_i,\;\; \theta^0=dt,\;\;
\theta^i=dx_i,\;\;\;(i=1,\cdots,n),
\end{align}
and $h:{\Bbb R}^{n+1}\to H$ is a smooth function, with
\begin{eqnarray}\label{eq:11}
H=\left\{\left(\begin{array}{lll}
b_0      &                    & \bigcirc \\[-2mm]
         & \!\!\!\ddots\!\!\! &          \\[-1mm]
\bigcirc &                    & b_n
\end{array}\right) \;|\;b_0,\cdots,b_n\in{\Bbb R}^*\right\}.
\end{eqnarray}

This is an overdetermined problem in coframes $\omega$ and
$\theta$; and, we proceed the overdetermined reduction method of
Proposition 9.12 of \cite{[7]}. The two coframes $\omega$ and
$\theta$ are connected by the ${\rm GL}(n+1,{\Bbb R})-$valued
function $A$; that is
\begin{eqnarray}\label{eq:12}
\omega=A\,\theta,\qquad A=\left(\begin{array}{cccc}
1      & 0      & \cdots & 0 \\
-f_1   & 1      & \cdots & 0 \\
\vdots & \vdots &        & \vdots\\
-f_n   & 0      & \cdots & 1 \\
\end{array}\right)
\end{eqnarray}
The entries of matrix $gAh^{-1}$ which are equal to
\begin{eqnarray}\label{eq:13}
gAh^{-1}=\left(\begin{array}{cccc}
a/b_0                                 & 0                  & \ldots & 0                  \\[2mm]
-(1/b_0)\sum_{i=1}^n a_{1i} f_i & a_{11}/b_1 & \ldots & a_{1n}/b_n \\[2mm]
\vdots                                &\vdots              &        & \vdots             \\[2mm]
-(1/b_0)\sum_{i=1}^n a_{ni} f_i & a_{n1}/b_1 & \ldots & a_{nn}/b_n
\end{array} \right)
\end{eqnarray}
are zeroth-order invariants, and we can normalize them by equating
$gAh^{-1}$ to
\begin{eqnarray}\label{eq:14}
\left(\begin{array}{cccccc}
1  & 0  & 0 & \cdots & 0 & 0 \\
-1 & 1  & 0 & \cdots & 0 & 0 \\
0  & -1 & 1 & \cdots & 0 & 0 \\
\vdots&\vdots&\vdots&&\vdots&\vdots \\
0  & 0 & 0 & \cdots  & -1 & 1\\
\end{array} \right).
\end{eqnarray}
This, leads to the $a=b_0$,
\begin{align}\label{eq:15}
&b_i=a_{i,i}=b_0/f_i,&i&=1,\cdots,n,\nonumber \\[-3mm] \\[-2mm]
&a_{i,i-1}=-b_0/f_{i-1},&i&=2,\cdots,n\nonumber
\end{align}
and other $a_{ij}$s are zero. Furthermore, we have a family of
coframes $h_0\,\theta$, with
\begin{align}\label{eq:16}
h_0=b_0\left(\begin{array}{cccc}
1      & 0     & \cdots & 0 \\
0      & 1/f_1 & \cdots & 0 \\
\vdots &\vdots &        & \vdots\\
0      & 0     & \cdots & 1/f_n
\end{array}\right).
\end{align}
A fixed coframe in this family is given by $b_0=1$ or
\begin{align}\label{eq:17}
\tilde{\theta}=\big(dt,dx_1/f_1,\cdots,dx_n/f_n\big)^t.
\end{align}
The elements of $H$ which have $\theta$ as the $\tilde{\theta}$
orbit are $h=b_0I_n$ ($I_n$ is the $n\times n-$identity matrix).
This shows that
\paragraph{Theorem.}
{\em The web geometry of system (\ref{eq:1}) under the
pseudo-group (\ref{eq:2}) is equivalence to a equivalence problem
with coframe
\begin{align}\label{eq:18}
\omega^0=dt,\qquad \omega^i=dx_i/f_i,\quad (i=1,\cdots,n),
\end{align}
and structure group $G=\left\{a\,I_n\; |\;a\in{\Bbb
R}^*\right\}$.}
\section{Absorption}
We lift the coframe (\ref{eq:18}) to $U\times G$ as
\begin{align}\label{eq:20}
\theta^0=a\,dt,\qquad \theta^i=a\,dx_i/f_i,\quad (i=1,\cdots,n),
\end{align}

If $n=1$, then the structure equations on $U\times G$ are
\begin{align}\label{eq:20.1}
d\theta^0=\alpha\wedge\theta^0,\qquad d\theta^1=\alpha\wedge\theta^1,
\end{align}
and $\alpha=da/a$ is the basic Maurer-Cartan form on $G$; There is
not any essential coefficients in (\ref{eq:20.1}). We must
prolonged the problem to $U\times G\cong {\Bbb R}^3$, and so on.
But, we can solve this problem by direct computations. Let
$dx/dt=f(x)$ and $dX/dT=F(X)$ are tow given equation, and
$(X,Y)=\Phi(t,x)=(\varphi_0(t),\varphi_1(x))$ be a diffeomorphism
which sent $dx/dt=f(x)$ to $dX/dT=F(X)$. Then
\begin{align}\label{eq:20.2}
F(X)=F(\varphi_1(x))=\frac{d\varphi_1(x)}{d\varphi_0(t)}=\frac{\varphi_1'(x)}{\varphi_0'(t)}\,\frac{dx}{dt}=\frac{\varphi_1'(x)}{\varphi_0'(t)}\,f(x).
\end{align}
Therefore $\Phi(t,x)=(t,F^{-1}(L(x)))$, where $L(x)$ is an
indefinite integral of $1/f(x)$. Thus, in this case
\paragraph{Theorem.}
{\em Any two first order homogeneous ODEs are web equivalence.}

\medskip Now, let $n\geq2$, then we had deduced structure equations
on $U\times G$ as
\begin{align}\label{eq:21}
d\theta^0 &= \alpha\wedge\theta^0, & d\theta^i &=
\alpha\wedge\theta^i+\sum_{j=1}^n\frac{\ell_{ij}}{a}\,\theta^i\wedge\theta^j,
&&(i=1,\cdots,n).
\end{align}
where
\begin{align}\label{eq:22}
\ell_{i,j}=f_j.\frac{\partial(\ln|f_i|)}{\partial x_j},\qquad
(i,j=1,\cdots,n,\;\;j\neq i).
\end{align}
and $\alpha=da/a$ is the basic Maurer-Cartan form on $G$. Let us
reduce the Maurer-Cartan form $\alpha$ back to the base manifold
$U$ by replacing them by general linear combinations of coframe
elements $\alpha \; \longmapsto \; \sum_{i=0}^nz_i\,\theta^i$.
This, leads to the system
\begin{align}\label{eq:24}
d\theta^0 &= \sum_{i=1}^nz_i\,\theta^0\wedge\theta^i, & d\theta^i
&=\sum_{j=1}^n\left(\frac{\ell_{ij}}{a}-z_j\right)\,\theta^i\wedge\theta^j,
&& (i=1,\cdots,n).
\end{align}
The combination $\ell_{ij}/a=(\ell_{ij}/a-z_j)-(z_j)$ is
invariant, and so will also contribute to the essential torsion.
We can normalize $\ell_{12}/a=1$ by setting $a=\ell_{1,2}$. This
normalization have the effect of eliminating all the group
parameters, and so, with just one loop through the equivalence
method, we have found an invariant coframe:
\begin{align}\label{eq:26}
\theta^0=dt,\qquad \theta^i=dx_i/(\ell_{12}\,f_i),\quad
(i=1,\cdots,n).
\end{align}
\paragraph{Theorem.}
{\em For $n\geq2$, the web geometry of system (\ref{eq:1}) under
the pseudo-group (\ref{eq:2}) is equivalence to
$\{e\}-$equivalence problem with coframe (\ref{eq:26}), with
structure equations
\begin{align}\label{eq:27}
d\theta^0 &=-\sum_{i=1}^n\frac{\partial\ell_{12}}{\partial x_i}\,\theta^0\wedge\theta^i,\nonumber\\[-3mm] \\[-2mm]
d\theta^i
&=\sum_{j=1}^n\left(\frac{\ell_{ij}}{\ell_{12}}-\frac{\partial\ell_{12}}{\partial
x_j}\right)\theta^i\wedge\theta^j,\qquad
(i=1,\cdots,n).\nonumber
\end{align}}

\medskip By the theorem 8.22 of \cite{[7]}, we conclude that
\paragraph{Theorem.}
{\em For $n\geq2$, the web symmetry group of a system (\ref{eq:1}) under
the pseudo-group (\ref{eq:2}) is a finite dimensional Lie group of dimension at most $n+1$.}

\end{document}